\documentclass[11pt,reqno]{amsart}
\usepackage{amssymb,amscd,amsbsy,mathrsfs}
\usepackage{mathdots}
\newcommand{\g}{\gamma}
\renewcommand{\d}{\delta}

\newcommand{\s}{\sigma}

\newcommand{\f}{\varphi}

\newcommand{\G}{\Gamma}
\newcommand{\D}{\Delta}

\newcommand{\cp}{{\mathscr P}}

\newcommand{\T}{{\Bbb T}}
\newcommand{\pp}{{\Bbb P}}
\newcommand{\dd}{{\Bbb D}}
\newcommand{\R}{{\Bbb R}}
\newcommand{\Z}{{\Bbb Z}}

\newcommand{\0}{{\boldsymbol{0}}}

\newcommand{\bs}{\boldsymbol}

\newcommand{\bS}{{\boldsymbol S}}

\newcommand{\rf}[1]{(\ref{#1})}

\newcommand{\df}{\stackrel{\mathrm{def}}{=}}

\newcommand{\supp}{\operatorname{supp}}

\newcommand{\rank}{\operatorname{rank}}
\newcommand{\const}{\operatorname{const}}

\newcommand{\eeq}{\end{equation}}
\newcommand{\beq}{\begin{equation}}
\newcommand{\bay}{\begin{eqnarray}}
\newcommand{\ba}{\begin{align*}}
\newcommand{\ea}{\end{align*}}
\newcommand{\ey}{\end{eqnarray}}
\newcommand{\bey}{\begin{eqnarray*}}
\newcommand{\eey}{\end{eqnarray*}}

\newcommand{\be}{\infty}

\newcommand{\bl}{\blacksquare}

\newcommand{\Pf}{{\bf Proof. }}

\newtheorem{thm}{\hspace{\parindent}Theorem}[section]

\newtheorem{lem}[thm]{\hspace{\parindent}Lemma}

\theoremstyle{remark}

\newtheorem*{rem*}{Remark}

\newcommand{\ri}{{\rm i}}

\newcommand\fM{\frak M}

\newcommand{\Bbbone}{{\rm{1\mathchoice{\kern-0.25em}{\kern-0.25em}{\kern-0.2em}{\kern-0.2em}I}}}

\newcommand\mB{\mathcal{B}}

\begin{document}

\numberwithin{equation}{section}

\numberwithin{equation}{section}

\title{Triangular projection on $\boldsymbol{\bS_p,~0<p<1}$ and related inequalities}
\author{A.B. Aleksandrov and V.V. Peller}
\thanks{The research on \S\:3--\S\:4  is supported by 
Russian Science Foundation [grant number 18-11-00053].
The research on \S\:5 is supported by 
Russian Science Foundation [grant number 20-61-46016].
The research on \S\:6 is supported by by a grant of the Government of the Russian Federation for the state support of scientific research, carried out under the supervision of leading scientists, agreement  075-15-2021-602.
The research on \S\:5 is supported by RFBR [grant number 20-01-00209]}.
\thanks{Corresponding author: V.V. Peller; email: peller@math.msu.edu}


\

\begin{abstract}
In this paper we study properties of the triangular projection $\cp_n$ on the space of $n\times n$ matrices. The projection $\cp_n$ annihilates the entries of an $n\times n$ matrix below the main diagonal and leaves the remaining entries unchanged. We estimate the $p$-norms of $\cp_n$ as an operator on the Schatten--von Neumann class $\bS_p$
for $0<p<1$. The main result of the paper shows that for $p\in(0,1)$, the $p$-norms of
$\cp_n$ on $\bS_p$
behave as $n\to\be$ as $n^{1/p-1}$. This solves a problem posed by B.S. Kashin.

Among other results of this paper we mention the result that describes the behaviour of 
the $\bS_p$-quasinorms of the $n\times n$ matrices whose entries above the diagonal are equal to 1 while the entries below the diagonal are equal to 0.
\end{abstract} 

\maketitle


\setcounter{section}{0}
\section{\bf Introduction}
\setcounter{equation}{0}
\label{In}

\

This paper continues studying Schur multipliers of Schatten--von Neumann classes 
$\bS_p$ for $p\in(0,1)$. We refer the reader to \cite{AP1} and \cite{AP2} for earlier results.

Recall that for infinite matrices $A=\{a_{jk}\}_{j,k\ge0}$ and 
$B=\{b_{jk}\}_{j,k\ge0}$ the {\it Schur--Hadamard product} $A\star B$ of $A$ and $B$ is, by definition,
the matrix
$$
A\star B=\{a_{jk}b_{jk}\}_{j,k\ge0}.
$$
In the same way one can define the Schur--Hadamard product of two finite matrices of the same size.

Let $p>0$. We say that a matrix $A=\{a_{jk}\}_{j,k\ge0}$ is a {\it Schur multiplier} of the Schatten--von Neumann class $\bS_p$ if
$$
B\in\bS_p\quad\Longrightarrow A\star B\in\bS_p.
$$
We use the notation  $\fM_p$ for the space of Schur multipliers of $\bS_p$. 

We refer the reader to \cite{GK1} and \cite{BS2} for properties of Schatten-von Neumann classes $\bS_p$. It is well known that for $p\ge1$, the class $\bS_p$ is a Banach space and $\|\cdot\|_{\bS_p}$ is a norm. For $p<1$, the class $\bS_p$ is a quasi-Banach space. More precisely, for $p<1$, $\bS_p$ is a
{\it $p$-Banach space} and $\|\cdot\|_{\bS_p}$ is a {\it $p$-norm}, i.e.,
$$
\|T+R\|_{\bS_p}^p\le\|T\|_{\bS_p}^p+\|R\|_{\bS_p}^p,\quad T,~R\in\bS_p,
$$
see, for example, \cite{Pe2}, App. 1, for the proof of this triangle inequality for $p<1$.

It follows from the closed graph theorem that if $A\in\fM_p$, then
$$
\|A\|_{\fM_p}\df\sup\{\|A\star B\|_{\bS_p}:~B\in\bS_p,~\|B\|_{\bS_p}\le1\}<\be.
$$
It is easy to see that $\fM_p$ is a Banach space for $p\ge1$ and it is a $p$-Banach space for $0<p<1$.

Obviously, all finite matrices are Schur multipliers of $\bS_p$ for all $p>0$ and we can define the norm $\|\cdot\|_{\fM_p}$ (the $p$-norm $\|\cdot\|_{\fM_p}$ if $p<1$) as in the case of infinite matrices.

We denote by $\mB(\bS_p)$ the space of bounded linear transformers on $\bS_p$ and 
for ${\mathcal T}\in\mB(\bS_p)$, we put
$$
\|{\mathcal T}\|_{\mB(\bS_p)}=\sup\big\{\|{\mathcal T}T\|_{\bS_p}:~T\in\bS_p,~\|T\|_{\bS_p}\le1\big\}.
$$
Clearly, $\|\cdot\|_{\mB(\bS_p)}$ is a $p$-norm for $p\in(0,1)$.

In this paper {\it we solve a problem by B.S. Kashin on the behavior of $\|\cp_n\|_{\mB(\bS_p)}$ for $p\in(0,1)$}, where $\cp_n$ is the upper triangular projection on the space of $n\times n$ matrices, i.e., for $A=\{a_{jk}\}_{1\le j\le n,1\le k\le n}$, the matrix 
$\cp_nA$ is defined by
$$
\cp_nA=\chi_n\star A,
$$
where the matrix $\chi_n\df\{(\chi_n)_{jk}\}_{1\le j\le n,1\le k\le n}$ is given by
$$
(\chi_n)_{jk}=\left\{\begin{array}{ll}1,&j\le k,\\0,&j>k.\end{array}\right.
$$

Let us start introduce the infinite Hankel matrix $\D_n\df\{(\D_n)_{jk}\}_{j,k\ge0}$ is defined by 
\bay
\label{beskDn}
(\D_n)_{jk}=\left\{\begin{array}{ll}1,&j+k<n,\\0,&j+k\ge n.\end{array}\right.
\ey
Clearly, $\|\chi_n\|_{\fM_p}=\|\D_n\|_{\fM_p}$.

The main result of the paper shows that for $0<p<1$, there exist positive numbers $c$ and $C$ such that
$$
cn^{1/p-1}\le\|\cp_n\|_{\mB(\bS_p)}=\|\chi_n\|_{\fM_p}=\|\D_n\|_{\fM_p}\le Cn^{1/p-1}.
$$
This solves a problem posed by B.S. Kashin.

Note that in the case when $1<p<\be$, the projections $\cp_n$, $n\ge1$, are uniformly bounded on 
$\bS_p$. This follows immediately from the well known fact that for $p\in(1,\be)$, the {\it triangular projection} $\cp$ onto the infinite upper triangular matrices is bounded, see \cite{GK2}. 
On the other hand,  in the case when $p=1$ the projection $\cp$ 
is unbounded, see \cite{GK2}.

Moreover, it is also well known that there are positive numbers $k$ and $K$ such that
\bay
\label{a_chto_esli_p=1?}
k \log(1+n)\le\|\chi_n\|_{\fM_1}\le K \log(1+n),
\ey
see Remark 6 in \S\;\ref{Osnova}. 

We give three different approaches to the problem of studying the behaviour of $\|\cp_n\|_{\mB(\bS_p)}$. The first approach is given in  \S\;\ref{Osnova}. It is based on Hankel operators of class $\bS_p$ and Hankel Schur multipliers of $\bS_p$.

The second approach is given in \S\;\ref{drugoi} while the third approach is given in \S\,\ref{eshchyo}.

In \S\;\ref{Sp-normy} we study the behaviour of the $p$-norms $\|\chi_n\|_{\bS_p}$, $p\le1$.
Note that this result will be used in \S\,\ref{eshchyo}.

A brief introduction in the Besov spaces $B_p^{1/p}$ will be given in \S\;\ref{Besovy}.
In \S\,\ref{Gankeli} we introduce Hankel matrices, certain Besov classes, we state a description of Hankel matrices of class $\bS_p$ in terms of Besov classes and also state upper estimates of $p$-Schur multiplier norms of Hankel matrices. Finally, in \S\;\ref{otenkitrigpolinomov} we give estimates of certain trigonometric polynomials that will be used later. 

\

\section{\bf The Besov classes $\bs{B_p^{1/p}}$}
\setcounter{equation}{0}
\label{Besovy}

\

For the reader's convenience we give here the definition of the Besov classes 
$\big(B_p^{1/p}\big)_+$ of functions in $B_p^{1/p}$ that are analytic in the unit disk $\dd$.

\medskip

{\bf Definition.} Let $v$ be an infinitely differentiable function on $(0,\be)$ such that
$v\ge\0$, $\supp v=\big[\frac12,2\big]$ and
$$
\sum_{j\ge0}v\big(2^{-j}x\big)=1\quad\mbox{for}\quad x\ge1.
$$
Let $0<p<\be$. For a function $\f$ analytic in $\dd$ we say that 
\bay
\label{Bp1/p}
\f\in\big(B_p^{1/p}\big)_+\quad\mbox{if}\quad
\sum_{n\ge0}2^n\|\f*V_n\|_{L^p}^p<\be,
\ey
where 
\bay
\label{Vn}
V_n(z)=\sum_{j>0}v\big(2^{-n}j\big)z^j
\quad\mbox{and}\quad
\f*V_n=\sum_{j>0}\widehat\f(j)\widehat V_n(j)z^j.
\ey
Among various equivalent norms on $\big(B_p^{1/p}\big)_+$ we can select the following one:
\bay
\label{normaBesova}
\|\f\|_{B_p^{1/p}}\df\left(\sum_{n\ge0}2^n\|\f*V_n\|_{L^p}^p\right)^{1/p}.
\ey

\

\section{\bf Hankel matrices and $L^p$ estimates for certain families of polynomials}
\setcounter{equation}{0}
\label{Gankeli}

\

In this section we recall certain facts on Hankel matrices (Hankel operators) that will be used in \S\;\ref{Osnova}. For a function $f$ analytic in $\dd$, we denote by $\G_\f$ the {\it Hankel matrix}
$$
\G_\f=\{\widehat\f(j+k)\}_{j,k\ge0},
$$
where $\widehat\f(j)$ is the $j$th Taylor coefficient of $\f$.

By Nehari's theorem, the matrix $\G_\f$ induces a bounded linear operator on the sequence space $\ell^2$ if and only if there exists a function $\psi$ in $L^\be(\T)$ such that their Fourier coefficients satisfy the equalities $\widehat\f(j)=\widehat\psi(j)$, $j\ge0$. Here $\widehat\psi(j)$ is the $j$th Fourier coefficient of $\psi$.

We will need the following criterion for the membership of $\G_\f$ in the Schatten--von Neumann class $\bS_p$: 

\medskip

{\it Let $0<p<\be$ and let $\f$ be a function analytic in $\dd$. Then
\bay
\label{GankeliBesovy}
\G_\f\in\bS_p\quad\Longleftrightarrow\quad\f\quad\mbox{belongs to the Besov class}\quad
B_p^{1/p}.
\ey}

This was established in
\cite{Pe1} for $p\ge1$. For $p<1$, this was proved in \cite{Pe2}, see also \cite{Pek} and \cite{S} for alternative proofs. We also refer the reader to the book \cite{Pe3}, see Ch.6, \S\,3, in which a detailed presentation of this material is given. 

The following result yields sharp estimates of $\|\G_\f\|_{\bS_p}$ for a certain class of polynomials $\f$ in terms of $\|\f\|_{L^p}$.

\begin{thm}
\label{spetsform}
Let $0<p<1$ and let $n$ be a positive integer. Suppose that $\f$ is a polynomial of the form
$$
\f(z)=\sum_{j=2^{n-1}+1}^{2^{n+1}-1}\widehat\f(j)z^j.
$$
Then
\bay
\label{s_dvukh_storon}
\d_p2^{(n+1)/p}\|\f\|_{L^p}\le\|\G_\f\|_{\bS_p}\le2^{(n+1)/p}\|\f\|_{L^p}
\ey
for a certain positive number $\d_p$.
\end{thm}

\Pf The upper estimate in \rf{s_dvukh_storon} is an immediate consequence of 
\rf{dlya_polinomov}. Let us proceed to the lower estimate. It follows from 
\rf{GankeliBesovy} and \rf{Bp1/p} that there exists a positive number $d_p$ such that
\bay
\label{f*Vk}
2^k\|\f*V_k\|_{L^p}^p\le d_p\|\G_\f\|_{\bS_p}^p,\quad k\ge1.
\ey
It follows from the definition of $V_k$ given in \rf{Vn} that for $n\ge2$,
$$
\f=\f*V_{n-1}+\f*V_n+\f*V_{n+1}.
$$
Thus, we can conclude from \rf{f*Vk} that
$$
\|\f\|_{L^p}^p\le\sum_{j=-1}^1\|\f*V_{n+j}\|^p_{L^p}\le3\cdot2^{n+1}d_p\|\G_\f\|_{\bS_p}^p.
$$
This implies the lower estimate for $\|\G_\f\|_{\bS_p}$ in \rf{s_dvukh_storon}. $\bl$

\medskip

{\bf Remark 1.} One can slightly modify the above reasoning and show that Theorem \ref {spetsform} can be extended to arbitrary $p>0$.

\medskip

{\bf Remark 2.} In the same way we can prove the following more general fact.
Suppose that $M>1$ and $k$ is a positive integer. Let $\f$ be an analytic polynomial of the form
$$
\f(z)=\sum_{k\le j<kM}\widehat\f(j)z^j.
$$
Then there are positive numbers $c_1$ and $c_2$ such that
$$
c_1(kM)^{(1/p-1)}\|\f\|^p_{L^p}\le\|\G_\f\|^p_{\bS_p}\le c_2(kM)^{(1/p-1)}\|\f\|^p_{L^p}.
$$

\medskip
Note that to prove the sufficiency of the condition $\f\in B_p^{1/p}$ 
for the membership of the Hankel operator $\G_\f$ in $\bS_p$
in the case $0<p\le1$, the following inequality was established in \cite{Pe1} and \cite{Pe2}, see also
\cite{Pe3}, Ch.6, \S\,3:
\bay
\label{dlya_polinomov}
\|\G_\f\|_{\bS_p}\le2^{1/p-1}m^{1/p}\|\f\|_{L^p}
\ey
for an arbitrary analytic polynomial $\f$ of degree at most $m-1$

\medskip

Let us proceed to Hankel Schur multipliers of $\bS_p$ for $p\in(0,1]$.
Let $0<p\le1$ and let $\f$ be an analytic polynomial of degree at most $m-1$. We need the following inequality:
\bay
\label{multnorm}
\|\G_\f\|_{\fM_p}\le(2m)^{1/p-1}\|\f\|_{L^p}.
\ey
It was established in \cite{Pe2}; the reader can also find the proof of this estimate in \cite{Pe3}, Ch.6, \S\,3.

\

\section{\bf  $\bs{L^p}$ estimates for certain trigonometric polynomials}
\setcounter{equation}{0}
\label{otenkitrigpolinomov}

\

We start this section with the following elementary remark.

\medskip

{\bf Remark 3.} It follows from the well-known inequality $|f(0)|\le\|f\|_{H^p}$, $p\in(0,\be]$,  that for an arbitrary trigonometric polynomial $f$ of the form
$f(z)=\sum\limits_{j=m}^n\widehat f(j)z^j$, the following inequalities hold:
\bay
\label{pervyi_i_poslednii}
|\widehat f(m)|\le\|f\|_{L^p}\quad\mbox{and}\quad |\widehat f(n)|\le\|f\|_{L^p}\quad\mbox{for all} \quad p\in(0,\be].
\ey

\medskip

Consider now the analytic Dirichlet kernel $D_n^+$ defined by
$$
D_n^+(z)=\sum_{j=0}^{n-1}z^j,\quad n\ge1.
$$

\begin{lem}
\label{Dirichlet}
Let $0<p<1$. There exists a positive number $K_p$ such that
$$
1\le\|D_n^+\|_{H^p}\le K_p
$$
for every $n\ge1$.
\end{lem}

\Pf By \rf{pervyi_i_poslednii}, $\|D_n^+\|_{H^p}\ge|D_n^+(0)|=1$.
On the other hand, it is easy to see that
$D_n^+(z)=\dfrac{1-z^n}{1-z}$. 
Hence,
$$
\|D_n^+\|_{H^p}^p\le\|(1-z)^{-1}\|_{H^p}^p+\|z^n(1-z)^{-1}\|_{H^p}^p=2\|(1-z)^{-1}\|_{H^p}^p.
$$
It remains to put $K_p=2\|(1-z)^{-1}\|_{H^p}^p$. $\bl$

\medskip

We also need an inequality for the $L^p$ norms of a family of trigonometric polynomials constructed from a smooth function on $\R$ with compact support. 

For an infinitely differentiable function $q$ on $\R$ with compact support, we define the family $Q_m$, $m\ge1$, of trigonometric polynomials defined by
\bay
\label{Ups}
Q_m(z)=\sum_{k\in\Z}q\left(\frac km\right)z^k,\quad m\ge1.
\ey
We need the following fact, which can be found in \cite{A}, Ch 4, \S\,2, see also \cite{Pe3}, Ch. 6, \S\,3:

\medskip

Suppose that $0<p\le1$ and $q$ is an infinitely differentiable function on $\R$ with compact support. Then there exists a positive number $C$ such that
\bay
\label{razduli}
\|Q_m\|_{L^p}\le Cm^{1-1/p},
\ey
where the trigonometric polynomial $Q_m$ is defined by \rf{Ups}.

\medskip

We are going to deduce from \rf{razduli} 
a lemma that shows to what extent the $L^p$-quasinorm of a trigonometric polynomial of a given degree can jump when the Riesz projection is applied. Recall that the Riesz projection 
$\pp_+$ and $\pp_-$ on the set of trigonometric polynomials are defined by
$$
\pp_+\psi(z)=\sum_{j\ge0} \widehat{\psi}(j)z^j\quad\mbox{and}
\quad\pp_-\psi(z)=\sum_{j<0} \widehat{\psi}(j)z^j.
$$

\begin{lem}
\label{rubim}
Let $0<p<1$ and let $q$ be an infinitely differentiable function on $\R$ such that 
$$
\supp q=[-1,1],\quad q(t)>0\quad\mbox{for}\quad t\in(-1,1)\quad\mbox{and}\quad q(0)=1.
$$ 
Consider the trigonometric polynomials  $Q_m$, $m\ge1$, of degree $m-1$ defined by {\em\rf{Ups}}. Then there is a positive number $d_p$ and a sequence of trigonometric polynomials $Q_m$, $m\ge1$, of degree $m-1$ such that
\bay
\label{rastyot}
\frac{\|\pp_+Q_m\|_{L^p}}{\|Q_m\|_{L^p}}\ge d_p m^{1/p-1}.
\ey
\end{lem}

\medskip 

\Pf Clearly, 
$\deg Q_m=m-1$ and $(\pp_+Q_m)(0)=1$. It is also easy to see that
$$
\|\pp_+Q_m\|_{L^p}\ge|(\pp_+Q_m)(0)|=1.
$$
On the other hand, by \rf{razduli}, $\|Q_m\|_{L^p}\le\const m^{1-1/p}$ which implies the result. $\bl$

\medskip

\medskip

{\bf Remark 4.} Actually, the inequality in \rf{rastyot} is optimal. Namely, if $0<p<1$, then there exists a positive number $C$ such that for an arbitrary trigonometric polynomial $Q$ of  degree $m-1$, $m\ge1$, the following inequality holds:
\bay
\label{Rieszp<1}
\|\pp_+Q\|_{L^p}\le C m^{1/p-1}\|Q\|_{L^p}.
\ey
Indeed, it is easy to see that it suffices to show that if $P$ is a trigonometric polynomial of the form
$$
P(z)=\sum_{j=2^{k-1}}^{2^{k-1}+2^k}\widehat P(j)z^j
\quad\mbox{and}\quad
R(z)=\sum_{j=2^{k-1}}^{2^k}\widehat P(j)z^j,
$$
then
$$
\|R\|_{L^p}\le\eta_p 2^{(1/p-1)k}\|P\|_{L^p}
$$
for a certain positive number $\eta_p$. Clearly,
$$
R=P*D^+_{2^k+1}\quad\mbox{and}\quad\G_R=\G_P\star\G_{D^+_{2^k+1}}.
$$ 
It follows from \rf{dlya_polinomov} that
$$
\|\G_R\|_{\bS_p}\le2^{(k+1)(1/p-1)}\|D^+_{2^k+1}\|_{L^p}\|\G_P\|_{\bS_p}
\le K_p2^{(k+1)(1/p-1)}\|\G_P\|_{\bS_p}
$$
by Lemma \ref{Dirichlet}. Inequality \rf{Rieszp<1} follows now from Theorem \ref{spetsform}.

Note that \rf{Rieszp<1} can be deduced from Theorem 1.1 of \cite{A}, Ch. 5.

%
%
%
%
%

%
%

\

\section{\bf The behaviour of $\bs{\|\D_n\|_{\bS_p}}$}
\setcounter{equation}{0}
\label{Sp-normy}

\

We find in this section the behaviour of the $\bS_p$-quasinorms of the Hankel matrices $\D_n$
for large values of $n$.

\begin{thm}
\label{DnSp}
Let $0<p\le1$. Then there exist positive numbers $c$ and $C$ such that 
$$
cn^{1/p}\le\|\D_n\|_{\bS_p}\le Cn^{1/p}.
$$
\end{thm}

It is an immediate consequence of \rf{GankeliBesovy} that Theorem \ref{DnSp} is equivalent to the following result.

\begin{lem}
\label{Dn+Besov}
Let $0<p\le1$. Then there exist positive numbers $c$ and $C$ such that
\bay
\label{Dir&Besov}
cn^{1/p}\le\|D^+_n\|_{B_p^{1/p}}\le Cn^{1/p}.
\ey
\end{lem}

{\bf Proof of Lemma \ref{Dn+Besov}.} First of all, it can be shown easily that it suffices to prove the desired inequalities in the case $n=2^k+1$. We have
$$
\big\|D^+_{2^k+1}\big\|_{B_p^{1/p}}^p=\sum_{j\ge0}2^j\big\|D^+_{2^k+1}*V_j\big\|_{L^p}^p
=\sum_{j=0}^{k-1}2^j\big\|V_j\big\|_{L^p}^p+2^k\big\|D^+_{2^k+1}*V_k\big\|_{L^p}^p,
$$
see \rf{normaBesova}. 

Let us first obtain the lower estimate for $\|D^+_n\|_{B_p^{1/p}}$. Clearly,
$$
D^+_{2^k+1}*V_k=\sum_{j=0}^{2^k}v\Big(\frac j{2^k}\Big)z^j,
$$
and so by \rf{pervyi_i_poslednii},
$$
2^k\big\|D^+_{2^k+1}*V_k\big\|_{L^p}^p\ge2^k|v(1)|^p=2^k
$$
which proves the lower estimate in \rf{Dir&Besov}.

Let us proceed now to the upper estimate in \rf{Dir&Besov}. First of all, it follows from 
\rf{razduli} that 
$$
2^j\big\|V_j\big\|_{L^p}^p\le\const2^j2^{j(p-1)}=\const2^{jp},
$$
and so
\bay
\label{pervye_k}
\sum_{j=0}^{k-1}2^j\big\|V_j\big\|_{L^p}^p\le\const\sum_{j=0}^{k-1}2^{jp}
\le\const2^{kp}.
\ey
Finally, it remains to estimate from above $2^k\big\|D^+_{2^k+1}*V_k\big\|_{L^p}^p$. By
\rf{s_dvukh_storon},
\begin{align*}
\big\|D^+_{2^k+1}*V_k\big\|_{L^p}
&\le\const2^{-(k+1)/p}\Big\|\G_{D^+_{2^k+1}*V_k}\Big\|_{\bS_p}\\[.2cm]
&=\const2^{-(k+1)/p}\Big\|\G_{D^+_{2^k+1}}\star\G_{V_k}\Big\|_{\bS_p}\\[.2cm]
&\le\const2^{(k+1)(1/p-1)}\|V_k\|_{L^p}
\Big\|\G_{D^+_{2^k+1}}\Big\|_{\bS_p}
\end{align*}
the last inequality being a consequence of \rf{multnorm}. By \rf{razduli},
$$
\|V_k\|_{L^p}\le\const2^{(k+1)(1-1/p)},
$$
and so by \rf{dlya_polinomov} and Theorem\ref{spetsform},
$$
\big\|D^+_{2^k+1}*V_k\big\|_{L^p}\le\const2^{-k/p}\Big\|\G_{D^+_{2^k+1}}\Big\|_{\bS_p}\le
\const\|D^+_{2^k+1}\|_{L^p}\le\const,
$$
the last inequality being a consequence of Lemma \ref{Dirichlet}. Thus,
$$
2^k\big\|D^+_{2^k+1}*V_k\big\|_{L^p}^p\le\const2^k.
$$
This together with \rf{pervye_k} proves the upper estimate in \rf{Dir&Besov}. $\bl$

\medskip

{\bf Remark 5.} It is easy to see that for $p\in(1,\be)$, there exist positive numbers $c$ and $C$ such that
$$
cn\le\|\D_n\|_{\bS_p}\le Cn.
$$

\

\section{\bf The main result}
\setcounter{equation}{0}
\label{Osnova}

\

In this section we prove the main result of the paper that solves a problem by B.S. Kashin.
It will be deduced from Theorem \ref{spetsform} and from Lemma \ref{rubim} below, which gives us an optimal estimate of the $L^p$ norms of trigonometric polynomials under the Riesz projection.

Recall that for $n\ge1$, $\|\chi_n\|_{\fM_p}=\|\D_n\|_{\fM_p}$,
where the infinite Hankel matrix $\D_n\df\{(\D_n)_{jk}\}_{j,k\ge0}$ is defined by \rf{beskDn}.

The main result of the paper is the following theorem:

\begin{thm}
\label{glavnaya}
Let $0<p<1$. 
There exist positive numbers $c_p$ and $C_p$ such that
$$
c_pn^{1/p-1}\le\|\chi_n\|_{\fM_p}=\|\D_n\|_{\fM_p}\le C_pn^{1/p-1}.
$$
\end{thm}


\medskip

{\bf Proof of Theorem \ref{glavnaya}.} Clearly, $\D_n$ is nothing but the Hankel matrix
$\G_{D_n^+}$. By \rf{multnorm}, 
$$
\|\D_n\|_{\fM_p}\le2^{1/p-1}n^{1/p-1}\|D_n^+\|_{L^p}\le C_pn^{1/p-1},
$$
where $C_p\df K_p2^{1/p-1}$.

Let us proceed to estimating $\|\D_n\|_{\fM_p}$ from below. Let $0<p<1$.
Clearly, it suffices to prove that 
$$
\|\chi_{2^k+1}\|_{\fM_p}=\|\D_{2^k+1}\|_{\fM_p}\ge\s_p2^{k(1/p-1)}
$$
for a certain positive number $\s_p$.

Consider the sequence of polynomials $Q_m$ constructed in Lemma \ref{rubim}.
Put
$$
P_k(z)=z^{2^k}Q_{2^{k-1}}(z),\quad k\ge1.
$$
Then $P_k$ is an analytic polynomial and 
$$
P_k(z)=\sum_{j=2^{k-1}}^{2^{k-1}+2^k}\widehat{P_k}(j)z^j.
$$
Put
$$
R_k(z)=\sum_{j=2^{k-1}}^{2^k}\widehat{P_k}(j)z^j.
$$
It is easy to see that
$
R_k=P_k*D^+_{2^k+1}.
$
Thus,
$$
\G_{R_k}=\G_{P_k}\star\D_{2^k+1}.
$$
It follows easily from \rf{rastyot} that
$$
\|R_k\|_{H^p}\ge\g_p2^{k(1/p-1)}\|P_k\|_{H^p}
$$
for a certain positive number $\g_p$. 

This together with Theorem \ref{spetsform} implies that
$$
\Big\|\D_{2^k+1}\star\G_{P_k}\Big\|_{\bS_p}
=\big\|\G_{R_k}\|_{\bS_p}\ge\rho_p2^{k(1/p-1)}\big\|\G_{P_k}\|_{\bS_p}
$$
for a certain positive number $\rho_p$. $\bl$

\medskip

{\bf Remark 6.} The above reasoning also allows one to prove the well known inequalities \rf{a_chto_esli_p=1?}. Indeed, to establish \rf{a_chto_esli_p=1?}, we can argue in the same way as in the case $p=1$, use the well known fact that the norms $\|D_n^+\|_{L^1}$
grow logarithmically and use instead of Lemma \ref{rubim} the well-known estimate
$$
\frac{\|\pp_+K_m\|_{L^1}}{\|K_m\|_{L^1}}\ge\const\log(1+m)
$$
for the Fej\'er kernel $K_m$. Also, we would like to mention that inequality \rf{multnorm} in the case $p=1$ is much easier than in the case $p<1$. Indeed, for $p=1$ it is well-known and follows from the implication
$$
\f\in L^1\quad\Longrightarrow\quad\G_\f\in\fM_1
$$
and
$$
\|\G_\f\|_{\fM_1}\le\|\f\|_{L^1},
$$
see \cite{Be}.

\

\setcounter{equation}{0}
\label{alternativa}
\section{\bf An alternative approach to the behaviour of  $\bs{\|\D_n\|_{\fM_p}}$}
\label{drugoi}

\

In this section we give an alternative approach to Theorem \ref{glavnaya}, i.e., we are going to give an alternative proof of the following inequalities:
\bay
\label{cpCp}
c_p n^{\frac1p-1}\le\|\chi_n\|_{\fM_p}\le C_p n^{\frac1p-1}
\ey
for certain positive numbers all $c_p$ and $C_p$ and for all positive integers $n$.
Clearly, it suffices to consider the case where $n=2^k$ for $k\in\Bbb N$.

We fix $p$ in $(0,1)$. In what all constants in inequalities can depend on $p$.

We need the following (apparently well-known estimate):
\bay
\label{snizu}
\|D_n^+\|_{B_p^{1/p}}\ge \const n^{\frac1p},\quad n\in{\Bbb N};
\ey
here the constant depend on $p$. Actually, it is also true that $\|D_n^+\|_{B_p^{1/p}}$ can be estimated from above in terms of $n^{\frac1p}$ but we do not need the upper estimate in this paper.


Applying \rf{pervyi_i_poslednii} for $f=D_{2^k+1}^+\star V_{2^k}$ and $n=2^k$, we find that
$$
\|D_{2^k+1}^+\star V_{2^k}\|_{H^p}\ge|\widehat{D_{2^k+1}^+\star V_{2^k}}(2^k)|=v(1)=1.
$$

It remains to observe that 
$$
\|D_{2^k+1}^+\|_{B_p^{1/p}}^p\ge\const\sum_{n\ge1}2^n\|D_{2^k+1}^+*V_n\|_{L^p}^p\ge\const2^k\|D_{2^k+1}^+*V_k\|_{L^p}^p
\ge\const 2^k.
$$
This certainly implies \rf{snizu}.

\medskip

{\bf An alternative proof of \rf{cpCp}.} Let us first estimate $\|\chi_n\|_{\fM_p}$ from below.
It follows from \rf{snizu} and  \rf{dlya_polinomov} that
$\big\|\G_{D_{2^k+1}^+}\big\|_{\bS_p}\ge\const2^{k/p}$, i.e., 
\bay
\label{chi2k}
\|\chi_{2^k}\|_{\bS_p}\ge\const2^{k/p}.
\ey

Denote by $\Bbbone_n$ the $n\times n$ matrix whose entries are identically equal to $1$.
Clearly, $\chi_{2^k}\star\Bbbone_{2^k}=\chi_{2^k}$.
Hence, 
$$
c_p2^{k/p}\le\|\chi_{2^k}\|_{\bS_p}\le\|\chi_{2^k}\|_{\fM_p}\|\Bbbone_{2^k}\|_{\bS_p}=2^k\|\chi_{2^k}\|_{\fM_p},
$$
and so by \rf{chi2k},
$$
2^k\|\chi_{2^k}\|_{\fM_p}\ge\const2^{k/p}
$$
which implies the left inequality in \rf{cpCp}.

We need the following equality which is a very special case of Theorem 3.2 in \cite{AP1}:
\bay
\label{ochen'chastnyi}
\left
\|\left(
\begin{array}{cc}
M & \0 \\
\0 & M \\
\end{array}
\right)\right\|_{\fM_p}=2^{\frac1p-1}\|M\|_{\fM_p}
\ey
for an arbitrary matrix $M$ and for $p\in(0,1)$.




 
\begin{lem} 
\label{corpm1}
Let $0<p<1$. Then 
$$
\|\chi_{2n}\|_{\fM_p}^p\le2^{1-p}\|\chi_n\|_{\fM_p}^p+1
$$
for any positive integer $n$.
\end{lem}

\Pf Clearly,
$$
\chi_{2n}=\left(
\begin{array}{cc}
\chi_n & \Bbbone_n \\
\0 & \chi_n \\
\end{array}
\right)
=\left(
\begin{array}{cc}
\chi_n & \0 \\
\0 & \chi_n \\
\end{array}
\right)
+
\left(
\begin{array}{cc}
\0 & \Bbbone_n \\
\0 & \0 \\
\end{array}
\right).
$$
It remains to apply \rf{ochen'chastnyi}. $\bl$

\medskip

Let us prove now the right inequality in \rf{cpCp}.
Applying Lemma \ref{corpm1} for $n=2^k$, we get
$$
\frac{2^{(p-1)(k+1)}\|\chi_{2^{k+1}}\|_{\fM_p}^p}{2^{(p-1)k}\|\chi_{2^k}\|_{\fM_p}^p}=\frac{\|\chi_{2^{k+1}}\|_{\fM_p}^p}{2^{1-p}
\|\chi_{2^k}\|_{\fM_p}^p}\le1+2^{p-1}\|\chi_{2^k}\|_{\fM_p}^{-p}
$$
for all positive integers $k$.
It follows that
$$
{2^{N(p-1)}\|\chi_{2^N}\|_{\fM_p}^p}=\prod_{k=0}^{N-1}\frac{2^{(p-1)(k+1)}\|\chi_{2^{k+1}}\|_{\fM_p}^p}{2^{(p-1)k}\|\chi_{2^k}\|_{\fM_p}^p}
\le\prod_{k=0}^\be\big(1+2^{p-1}\|\chi_{2^k}\|_{\fM_p}^{-p}\big).
$$
It remains to observe that the infinite product converges by the left inequality in \rf{cpCp} which has already been proved above. $\bl$

\

\setcounter{equation}{0}
\label{eshchyo}
\section{\bf Yet another approach}
\label{drugoi}

\

The purpose of this section is to give yet another approach to describe the behaviour
of the $p$-norms $\|\cp_n\|_{\mB(\bS_p)}$ for large values of $n$. We need the following well known result: 

\medskip

{\it The operator of triangular projection $\cp$has weak type (1,1), i.e.,
\bay
\label{slabo!}
T\in\bS_1\qquad\Longrightarrow\qquad
s_j(\cp T)\le\const\|T\|_{\bS_1}(1+j)^{-1},\quad j\ge0,
\ey
where $\{s_j(\cp T)\}_{j\ge0}$ is the sequence of singular values of the operator $\cp T$.}

\medskip

Indeed, inequality \rf{slabo!} is an easy consequence of the following well known fact on Volterra operators with trace class imaginary parts:

\medskip

{\it Let $A=G+\ri H$ be a Volterra operator (i.e., $A$ is compact and its spectrum consists $\s(A)$ satisfies $\s(A)=\{0\}$). Then
\bay
\label{Volterra}
H\in\bS_1\qquad\Longrightarrow\qquad
s_j(G)\le\const\|H\|_{\bS_1}(1+j)^{-1},\quad j\ge0.
\ey
}

\medskip

This result can be found \cite{GK2}, Ch. 7, see also \cite{BS1}.

Let us explain for the sake of convenience how to deduce \rf{slabo!} from \rf{Volterra}.
Clearly, it suffices to establish \rf{slabo!} for self-adjoint operators $T$ with zero diagonal matrix entries, i.e., $t_{jj}=0$ for $j\ge0$, where $\{t_{jk}\}_{j,k\ge0}$ is the
matrix of $T$. Since $\cp$ is a bounded transformer on $\bS_p$ for $p\in(1,\be)$, it follows that
$\cp T\in\bS_p$ for all $p>1$. It is easy to verify that $\cp T$ is a Volterra operator and the imaginary part of $\cp(\ri T)$ is equal to
$$
\frac12\big(\cp(\ri T)-(\cp(\ri T))^*\big)=\frac\ri2\big(\cp T+(\cp T)^*\big)=\frac\ri2 T\in\bS_1.
$$
Then by \rf{Volterra}, \rf{slabo!} holds.

Let us now proceed to another alternative proof of Theorem \ref{glavnaya}.

\medskip

{\bf Yet another proof of Theorem \ref{glavnaya}.} Let us first show that for $p\le1$,
\bay
\label{verkhn1p-1}
\|\cp_n\|_{\mB(\bS_p)}\le\const n^{1/p-1}.
\ey
Suppose that $T$ is a rank one $n\times n$ matrix. Then the operator norm of $T$ coincides with $\|T\|_{\bS_p}$ for all $p>0$. Then by \rf{slabo!},
\begin{align}
\label{ranga1}
\|\cp_nT\|_{\bS_{p}}&=\left(\sum_{j=0}^{n-1}\big(s_j(\cp_nT)\big)^p
\right)^{1/p}\nonumber\\[.2cm]
&\le\const\left(\sum_{j=0}^{n-1}(1+j)^{-p}\right)^{1/p}\|T\|_{\bS_p}\le\const n^{1/p-1}\|T\|_{\bS_p}.
\end{align}
Suppose now that $T$ is an arbitrary $n\times n$ matrix. Then $T$ can be represented as
$
T=\sum_{j=1}^nT_j,
$
where $\rank T_j=1$ for $1\le j\le n$ and 
$$
\|T\|_{\bS_p}^p=\sum_{j=1}^n\|T_j\|_{\bS_p}^p.
$$
Hence, by \rf{ranga1},
\begin{align*}
\|\cp_nT\|^p_{\bS_{p}}&=\left\|\cp_n\left(\sum_{j=1}^nT_j\right)\right\|^p_{\bS_{p}}
\le\sum_{j=1}^n\|\cp T_j\|^p_{\bS_p}\\[.2cm]
&\le\const n^{1/p-1}\sum_{j=1}^n\|T_j\|^p_{\bS_p}=\const n^{1/p-1}\|T\|_{\bS_p}^p
\end{align*}
which proves \rf{verkhn1p-1}.

To prove that $\|\cp_n\|_{\mB(\bS_p)}\le\const n^{1/p-1}$, we apply the transformer
$\cp_n$ to the matrix $\Bbbone_n$ defined in \S\;\ref{alternativa}. Clearly, 
$\|\Bbbone_n\|_{\bS_p}=n$. On the other hand, by Theorem \ref{DnSp},
$$
\|\cp\Bbbone_n\|_{\bS_p}=\|\D_n\|_{\bS_p}\ge\const n^{1/p},
$$
and so
$$
\|\cp\Bbbone_n\|_{\bS_p}\ge\const n^{1/p-1}\|\Bbbone_n\|_{\bS_p}
$$
which proves the result. $\bl$

\

\
 
 \begin{footnotesize}
 
\noindent
\begin{tabular}{p{7cm}p{15cm}}
A.B. Aleksandrov & V.V. Peller \\
St.Petersburg Department & Department of Mathematics\\
Steklov Institute of Mathematics  & and Computer Sciences\\
Fontanka 27, 191023 St.Petersburg & St.Petersburg State University\\
Russia&Universitetskaya nab., 7/9,\\
email: alex@pdmi.ras.ru&199034 St.Petersburg, Russia\\
\\
&St.Petersburg Department\\
&Steklov Institute of Mathematics\\
&Russian Academy of Sciences\\
&Fontanka 27, 191023 St.Petersburg\\
&Russia\\
\\
&Department of Mathematics\\
&Michigan State University\\
&East Lansing, Michigan 48824\\
&USA\\

& email: peller@math.msu.edu
\end{tabular}

\end{footnotesize}

\end{document}